\documentclass[12pt]{amsart}

\newtheorem{theorem}{Theorem}[section]

\newtheorem{lemma}[theorem]{Lemma}
\theoremstyle{proposition}

\newtheorem{corollary}[theorem]{Corollary}

\newtheorem{definition}[theorem]{Definition}
\newtheorem{example}[theorem]{Example}

\newtheorem{question}[theorem]{Question}
\theoremstyle{remark}

\newtheorem{remark}[theorem]{Remark}

\numberwithin{equation}{section}

\usepackage{times}
\usepackage{enumerate}
\usepackage{graphicx,adjustbox}
\usepackage{hyperref}
\usepackage{amsmath,amssymb}
\usepackage{amscd}
\usepackage{graphicx}
\usepackage[all]{xy}
\usepackage{mathtools}

\title[Projective Closure of Affine Monomial Curves II]
{Projective Closure of Affine Monomial Curves II}
\author{
Joydip Saha
\and
Indranath Sengupta
\and
Pranjal Srivastava
}
\date{}
\address{\small \rm  Stat Math Unit, Indian Statistical Institute, Kolkata, West-Bengal,700108, INDIA.} 
\email{saha.joydip56@gmail.com}
\thanks{The first author thanks NBHM, Government of India for post-doc fellowship at ISI kolkata}

\address{\small \rm  Discipline of Mathematics, IIT Gandhinagar, Palaj, Gandhinagar, 
Gujarat 382355, INDIA.}
\email{indranathsg@iitgn.ac.in}
\thanks{The second author is the corresponding author.}

\address{\small \rm  Discipline of Mathematics, IIT Gandhinagar, Palaj, Gandhinagar, 
Gujarat 382355, INDIA.}
\email{pranjal.srivastava@iitgn.ac.in}
\date{}



\begin{document}

\begin{abstract}
In this paper our aim is twofold. First, we introduce the notion of star gluing of numerical 
semigroups and show that arithmetically Cohen-Macaulay and Gorenstein 
properties of the projective closure are preserved under this gluing operation. 
We then give a condition on Gr\"{o}bner basis 
of the defining ideal of an affine monomial curve which ensures that the 
Betti sequence of the affine curve is the same as the Betti sequence of its 
projective closure. We also study the effect of simple gluing on Betti 
sequences of the projective closure. Finally, we construct some numerical 
semigroups, using a gluing technique, such that the Cohen-Macaulay type 
of corresponding affine curve and its projective closure are both $n$.
\end{abstract}

\maketitle

\section{Introduction}
In general it is difficult to find a projective closure of an affine 
monomial curve, which is arithmetically Cohen-Macualay. For example, 
one such class can be found in \cite{ACM}, where Patil et.al. proved 
that the projective closure of the affine monomial curve defined by an 
arithmetic sequence is arithmetically Cohen-Macaulay. It is well-known 
that if the semigroup ring of a numerical semigroup is 
Cohen-Macaulay (respectively Gorenstein) then it is not always 
true that its projective closure is arithmetically Cohen-Macaulay 
(respectively Gorenstein). This makes us ask the first 
question, which is one of the main questions studied in this paper:
\medskip

\noindent\textbf{Question \ref{Question}} Suppose the projective closures of two affine monomial curves are arithmetically Cohen-Macaulay (respectively Gorenstein). Which condition(s) on gluing does preserve the arithmetically Cohen-Macaulay 
(respectively Gorenstein) property of the projective closure of the monomial curve obtained by gluing of these two monomial curves?
\medskip

An answer to the above Question would help us create 
a large family of affine monomial curves, with arithmetically 
Cohen-Macaulay (respectively Gorenstein) projective closure. 
By an abuse of terminology, we often say that a numerical 
semigroup is Cohen-Macaulay (respectively Gorenstein), 
which only means that the associated semigroup ring is 
Cohen-Macaulay (respectively Gorenstein). 
Similarly, we often use the terminology  "projective closure of a numerical semigroup" 
instead of the longer phrase "the projective closure of the 
affine monomial curve defined by the numerical semigroup". 
\medskip

Rosales \cite{JCR} introduced 
the concept of gluing of numerical semigroups, which was 
motivated by Delorme's work on 
complete intersection numerical semigroups. In \cite{JCR}, 
Rosales has given the minimal generating set of the defining 
ideal of the affine monomial curve associated to the gluing 
of numerical semigroups. Later, several 
mathematicians used this technique to create many examples 
of set-theoretic complete intersection and ideal-theoretic 
complete intersection of affine and projective varieties. 
Feza Arslan et al. \cite{GHM} used \textit{nice gluing} to 
give infinitely many families of $1$-dimensional local rings with non-Cohen 
Macaulay tangent cone and non-decreasing Hilbert function. 
\medskip

We introduce a new notion of gluing of two numerical semigroups and 
call it the \textit{star gluing}. Our aim is to understand, under what 
condition, are the properties of arithmetically Cohen-Macaulay and 
Gorenstein retained for the projective closure of an affine monomial 
curve defined by the star gluing of two numerical semigroups. We discuss 
this in section 3 after a brief discussion on terminologies and notations  
in section 2. 
\medskip

Our next question is about finding conditions on a numerical semigroup 
such that the Betti sequences of the affine monomial curve and its 
projective closure are the same. We ask the following question:
\medskip

\noindent\textbf{Question 2.} When can we say that the Betti numbers of the 
affine monomial curve and its projective closure are the same?
\medskip

In section 4, we will find a condition on Gr\"{o}bner basis of the defining ideal of the 
affine monomial curve such that the Betti numbers of the projective closure is the 
same as those of the affine monomial curve. An important class of monomial curves 
which satisfy this property is the one defined by arithmetic sequences. 
In section 5 we discuss  the Betti numbers of the projective closures of numerical 
semigroups obtained by the simple gluing of numerical semigroup. Section 6 deals with 
the construction of numerical semigroup, independent of the  embedding dimension, such 
that the corresponding affine curve and its projective closure both are Gorenstein.

\section{Preliminary}
Let $r\geq 3$ and $\mathbf{\underline n} = (n_{1}, \ldots, n_{r})$ be a sequence of 
$r$ distinct positive integers with $\gcd(\mathbf{\underline n})=1$. Let us assume 
that the numbers $n_{1}, \ldots, n_{r}$ generate the numerical semigroup 
$\Gamma(n_1,\ldots, n_r)$ in the set of all non-negative integers $\mathbb{N}$:
$$\Gamma(n_1,\ldots, n_r) = \lbrace\sum_{j=1}^{r}z_{j}n_{j}\mid z_{j'}s \text{ are non-negative \, integers}\rbrace$$ 
minimally, that is, if $n_i=\sum_{j=1}^{r}z_{j}n_{j}$ for some non-negative 
integers $z_{j}$, then $z_{j}=0$ for all $j\neq i$ and $z_{i}=1$. Let 
$\eta:k[x_1,\,\ldots,\, x_r]\rightarrow k[t]$ be the mapping defined by 
$\eta(x_i)=t^{n_i},\,1\leq i\leq r$, $\frak{p}(n_1,\ldots, n_r) = \ker (\eta)$ (called the defining ideal of $\Gamma(n_1,\ldots, n_r)$ ) and $\beta_{i}(\frak{p}(n_1,\ldots, n_r))$ denote the $i$-th Betti number of the ideal 
$\frak{p}(n_1,\ldots, n_r)$. Therefore, $\beta_{1}(\frak{p}(n_1,\ldots, n_r))$ denotes the 
minimal number of generators $\frak{p}(n_1,\ldots, n_r)$. We also denote the ring $k[x_{1},\ldots,x_{r}]/\frak{p}(n_1,\ldots, n_r)$ by $k[\Gamma(n_1,\ldots, n_r )]$, which is called the semigroup ring w.r.t. the semigroup $\Gamma(n_1,\ldots, n_r )$. It is known that $\frak{p}(\Gamma)$ is generated by binomials $x^{a}-x^{b}$ where $a$ and $b$ are $r$-tuples of non-negative integers with $\eta(x^{a})=\eta(x^{b})$. The set $\{x_{i} \vert a_{i}+b_{i} \neq 0 \}$ of all variables which appear in $f= x^{a}-x^{b}$, is called the support of $f$ and is denoted by $\mathrm{supp}(f)$.
\medskip

Let us assume that $n_{r}>n_{i}$ for all $i<r$, for the sequence $\mathbf{\underline n} = (n_{1}, \ldots, n_{r})$ and $n_{0}=0$. We define the semigroup $\overline{\Gamma(n_{1},\ldots n_{r})}\subset \mathbb{N}^{2}$, ( $\mathbb{N}$ is the set of all non negetive integers) generated by $\{(n_{i},n_{r}-n_{i})\mid 0\leq i\leq r\}$. Let us denote $\overline{\mathfrak{p}(n_{1},\ldots n_{r})}$ be the kernel of $k$-algebra map $\eta^{h}:k[x_{0},\ldots,x_{r}]\longrightarrow k[s,t]$, $\eta^{h}(x_{i})=t^{n_{i}}s^{n_{r}-n_{i}}, 0\leq r$. Then homogenization of the ideal $\frak{p}(n_1,\ldots, n_r) $ w.r.t. the variable $x_{0}$ is $\overline{\mathfrak{p}(n_{1},\ldots n_{r})}$. Thus the projective curve $\{[(a^{n_{r}}:a^{n_{r}-n_{1}}b^{n_{1}}:\cdots:b^{n_{r}})]\in\mathbb{P}^{n}_{k}\mid a,b\in k\}$ is the projective closure of the affine curve $C(n_{1},\ldots,n_{r}):=\{(b^{n_{1}},\ldots b^{n_{r}})\in \mathbb{A}^{n}_{k}\mid b\in k\}$ and we denote it by $\overline{C(n_{1},\ldots,n_{r})}$ (or by $\overline{C(\Gamma)}$, where $\Gamma=\Gamma(n_{1},\ldots n_{r})$) and the ring $k[x_{0},\ldots,x_{r}]/\overline{\frak{p}(n_1,\ldots, n_r)}$ by $k[\overline{\Gamma(n_1,\ldots, n_r )}]$. The Projective curve $\overline{C(n_{1},\ldots,n_{r})}$ is arithmetically Cohen-Macaulay if the vanishing ideal $\overline{\mathfrak{p}(n_{1},\ldots n_{r})}$ is a Cohen-Macaulay ideal. 
\bigskip
						
\section{Star gluing of numerical semigroups and its projective closures}
Let $\overline{\Gamma(n_1,\ldots, n_r )}=\langle(n_{1},n_{r}-n_{1}),\dots,(n_{r},0),(0,n_{r})\rangle$. 
Let us introduce two numerical semigroups $\Gamma_{1}$, $\Gamma_{2}$ as follows: 
$\Gamma_{1}=\langle n_{1},n_{2},\dots,n_{r} \rangle$ and 
$ \Gamma_{2}=\langle n_{r}-n_{1},\dots,n_{r}-n_{r-1},n_{r}\rangle$.

\begin{definition}\label{good}\cite{CN}{\rm 
Let $B_{i}$ be the Ap\'{e}ry set of $\Gamma_{i}(i=1,2)$ with respect to $n_{r}$.
\begin{itemize}
\item[(a)] We call Ap\'{e}ry set with respect to $n_{d}$ of $\overline{\Gamma(n_1,\ldots, n_r )}$ the generating set $$B=\{b_{0},b_{1},\dots,b_{n_{d}}\}$$ 
where $b_{0}=(0,n_{d}),b_{n_{d}}=(n_{d},0)$ and $b_{i}=(\upsilon_{i},\mu_{i})$ are such that:
\begin{itemize}
\item[(i)] $\{ \upsilon_{1},\dots,\upsilon_{n_{k-1}},n_{k}\}=B_{1}$.
\item[(ii)] $ \mu_{i}$ is the least element of $\overline{\Gamma(n_1,\ldots, n_r )}_{2}$ such that $(\upsilon_{i},\mu_{i})\in\overline{\Gamma(n_1,\ldots, n_r )}(1 \leq i \leq n_{k-1})$.
\end{itemize}

\item[(b)] The Ap\'{e}ry set $B$ is said to be good if $\{\upsilon_{1},\dots, \upsilon_{n_{k-1}},n_{k}\}=B_{2}$.
\end{itemize}
}
\end{definition}

\begin{lemma}\label{GD} The following conditions are equivalent:

\begin{itemize}
\item[(1)] $k[\overline{\Gamma(n_{1},\ldots,n_{r})}]$ is  arithmetically Cohen-Macaulay.
\item[(2)] $t^{n_{r}},s^{n_{r}}$ is a regular sequence.
\item[(3)] $((n_{r},0)+\overline{\Gamma(n_1,\ldots, n_r )})\cap ((0,n_{r})+\overline{\Gamma(n_1,\ldots, n_r )})=(n_{r},n_{r})+\overline{\Gamma(n_1,\ldots, n_r )}$.
\item[(4)] The Ap\'{e}ry Set $B$ of $ \overline{\Gamma(n_1,\ldots, n_r )}$ is good.
\end{itemize}

\end{lemma}
\proof Follows from the Lemma 4.3 and the Theorem 4.6 in \cite{CN}. \qed
\medskip
 
\begin{definition}(\cite{GHM}){\rm 
Let $\Gamma_{1} = \Gamma(m_{1},\dots, m_{l})$ and 
$\Gamma_{2} = \Gamma(n_{1},\dots, n_{k})$ be two numerical semigroups with $m_{1} < \cdots < m_{l}$ and $n_{1} < \cdots < n_{k}$. Let $p =b_{1}m_{1} +\cdots +b_{l}m_{l} \in \Gamma_{1}$ and $q = a_{1}n_{1} +\cdots +a_{k}n_{k} \in\Gamma_{2}$ be two positive integers satisfying
$\gcd(p, q) = 1 $ with $p \notin \{m_{1}, \dots ,m_{l} \}$, $q \notin \{n_{1}, \ldots ,n_{l}\}$ and $\{qm_{1},\ldots , qm_{l}\}\cap
\{pn_{1},\ldots , pn_{k}\} = \phi$. The numerical semigroup $\Gamma_{1}\#_{p,q} \Gamma_{2}=\langle qm_{1},\ldots , qm_{l}, pn_{1},\ldots, pn_{k}\rangle$ is
called a \textit{gluing} of the semigroups $\Gamma_{1}$ and $\Gamma_{2}$.
}
\end{definition}
\medskip

\begin{definition}{\rm 
The numerical semigroup $\Gamma_{1}\#_{p,q} \Gamma_{2}$,  
obtained by gluing of  $\Gamma_{1} = \Gamma(m_{1},\dots, m_{l})$ and $\Gamma_{2} = \Gamma(n_{1},\dots, n_{k})$ with respect to the positive integers $p$ and $q$ is 
called a \textit{star gluing} if $p = b_{l}m_{l} \in \Gamma_{1}$ and $q = a_{1}n_{1}+a_{2}n_{2}+\dots+a_{k}n_{k} \in \Gamma_{2}$, 
with $ a_{1}+a_{2}+\dots+a_{k} \leq b_{l}.$ We use the notation 
$\Gamma_{1}\star_{p,q} \Gamma_{2}$ to denote the star gluing 
of $\Gamma_{1}$ and $\Gamma_{2}$ with respect to the positive integers 
$p$ and $q$.
}
\end{definition}
\medskip

\begin{remark} It is important to find the largest integer from the generator of the numerical semigroup  $\Gamma_{1}\star_{p,q} \Gamma_{2}$ obtained by star gluing. The conditions $a_{1}+a_{2}+\dots+a_{k} \leq b_{l}$ with $m_{1}<\dots < m_{l},n_{1}<\dots <n_{k}$, $\gcd(p,q)=1$ and 
$\{qm_{1}, \dots, qm_{l}\}\cap \{ pn_{1}, \dots, pn_{k}\}=\phi$ imply that $$pn_{k}=(b_{l}m_{l})n_{k} \geq (a_{1}+a_{2}+\dots+a_{k})m_{l}n_{k}\geq (a_{1}n_{1}+a_{2}n_{2}+\dots+a_{k}n_{k})m_{l}=qm_{l}$$
and $pn_{k}$ is the largest integer among the generators of numerical semigroup  $\Gamma_{1}\star_{p,q} \Gamma_{2}$.
\end{remark}
\medskip

Suppose the projective closures of two monomial curves $\overline{C(m_{1},\dots, m_{l})}$ and $\overline{C(n_{1},\dots, n_{k})}$ are arithmetically Cohen-Macaulay, then the projective monomial curve $\overline{C(qm_{1},\ldots, qm_{l}, pn_{1},\ldots , pn_{k})}$ obtained by gluing of these two numerical semigroup may not be necessarily  arithmetically Cohen-Macaulay. We give the following example.
\medskip

\example {\rm We consider monomial curves $\overline{C(3,5)}$ and  $\overline{C(7,12)}$. They are Cohen-Macaulay (Gorenstein), but their gluing with respect to elements 
$p=8, q=19$, $\overline{C((57,95,56,96))}$ is neither arithmetically 
Cohen-Macaulay nor Gorenstein (verified with the help of MACAULAY 2 \cite{Macaulay}).
}
\medskip

\noindent This motivates us to ask the following question:
\medskip

\question\label{Question} {\rm Suppose the projective closures of two affine monomial curves are arithmetically Cohen-Macaulay (respectively Gorenstein). Which condition(s) on gluing does preserve the arithmetically Cohen-Macaulay 
(respectively Gorenstein) property of the projective closure of the monomial curve obtained by gluing of these two monomial curves?
}
\medskip

\noindent The following Lemma will be useful in proving our main theorems.
\medskip

\begin{lemma}\label{DIG} If $\Gamma$ is obtained by gluing $\Gamma_{1}$ and $ \Gamma_{2}$, if the defining ideals $\mathfrak{p}(\Gamma_{1}) \subset k[x_{1},\ldots, x_{l}]$ and $\mathfrak{p}(\Gamma_{2}) \subset k[y_{1},\ldots, y_{k}]$ are generated by the sets $G_{1} = \{f_{1},\ldots, f_{s}\} $ and $G_{2} = \{g_{1},\ldots, g_{t}\}$ respectively,
then the defining ideal $\frak{p}(\Gamma)\subset k[x_{1},\ldots, x_{l}, y_{1},\ldots, y_{k}]$ is generated by the set
$G = G_{1}\cup G_{2}\cup \{\rho\}$.
\end{lemma}
\proof: See Lemma 2.2 in \cite{SMR}. \qed
\medskip

\begin{lemma}\label{GG}
Let $\Gamma_{1} = \Gamma(m_{1},\dots, m_{l})$ and $\Gamma_{2} = \Gamma(n_{1},\dots, n_{k})$ be two numerical semigroups with $m_{1} < \cdots < m_{l}$ and $n_{1} < \cdots < n_{k}$. Let 
$p =b_{1}m_{1} +\cdots +b_{l}m_{l} \in \Gamma_{1}$ and 
$q = a_{1}n_{1} +\cdots +a_{k}n_{k} \in\Gamma_{2}$ be 
two positive integers satisfying 
$\gcd(p, q) = 1 $, such that $\Gamma=\Gamma_{1}\star_{p,q} \Gamma_{2}$ is the star gluing of $\Gamma_{1}$ 
and $\Gamma_{2}$. Suppose that $G_{1}$ and $G_{2}$ 
are Gr\"{o}bner bases of $\frak{p}(\Gamma_{1})$ and 
$\frak{p}(\Gamma_{2})$, with respect to the degree 
reverse lexicographic ordering induced by 
$x_{1}>\dots>x_{l}$ and $y_{1}>\dots>y_{k}$ 
respectively. Then 
$G_{1}^{h} \cup G_{2}^{h} \cup \{x^{b_{l}}_{l}-x_{0}^{b_{l}-(a_{1}+\cdots+a_{k})}y^{a_{1}}_{1}y^{a_{2}}_{2}\dots y^{a_{k}}_{k}\}$ is a Gr\"{o}bner bases of 
$\overline{\frak{p}(\Gamma)}$  with respect to the degree 
reverse lexicographic ordering 
$x_{1} > \dots>x_{l} > y_{1} > \dots > y_{k} > x_{0}$.
\end{lemma}

\proof We denote the $s$-polynomial of the polynomials $f$ and $g$ by $S(f, g)$. Let $G_{1} = \{f_{1},\ldots, f_{s}\}$ be a reduced
Gr\"{o}bner basis of the ideal $\mathfrak{p}(\Gamma_{1}) \subset k[x_{1},\ldots, x_{l}]$ with respect to the degree
reverse lexicographical ordering induced by $x_{1}> \cdots> x_{l}$ and $ G_{2} = \{g_{1}, \ldots , g_{t}\}$ be
a reduced Gr\"{o}bner basis of the ideal $ \mathfrak{p}(\Gamma_{2}) \subset k[y_{1}, \ldots , y_{k}]$ with respect to the degree
reverse lexicographical ordering induced by $y_{1}>\ldots > y_{k}$.
We will show that Gr\"{o}bner basis  of $\frak{p}(\Gamma)$ is  $$G=\{f_{1}, \dots, f_{s},g_{1}, \dots , g_{t},x_{l}^{b_{l}}-y^{a_{1}}_{1}y^{a_{2}}_{2}\dots y^{a_{k}}_{k} \}$$ with respect to degree
reverse lexicographical ordering $x_{1} >\cdots>x_{l}>y_{1}>\cdots>y_{k}$. By Lemma \ref{DIG}, defining ideal $\mathfrak{p}(\Gamma)$ of the affine curve $C(\Gamma)$ obtained by gluing is generated by the set $$G=\{f_{1}, \dots, f_{s},g_{1}, \dots , g_{t},x^{b_{l}}_{l}-y^{a_{1}}_{1}y^{a_{2}}_{2}\dots y^{a_{k}}_{k} \}.$$ 
Let $NF( -|G)$ be a weak normal form of a polynomial $(-)$ with respect to $G$ (definition 1.6.5 in \cite{BC}). If $ f,g \in k[x_{0},x_{1},\ldots,x_{l},y_{1},\ldots,y_{k}]$ be polynomial such that $$ \mathrm{lcm}(LM(f),LM(g))=LM(f)LM(g)$$ then $NF(S(f,g)|G)=0$. 
If we take the degree reverse lexicographic ordering $x_{1}>\cdots>x_{l}>y_{1}>\cdots>y_{k}$ then leading monomial of any element of $ G_{1}$ is product of $ x_{i}^{\alpha_{i}}$  and $G_{2}$ is product of $ y_{j}^{\beta_{j}}$ for some $i,j$ and $LM(x^{b_{l}}_{l}-y^{a_{1}}_{1}y^{a_{2}}_{2}\dots y^{a_{k}}_{k})=x_{l}^{b_{l}}$, so $NF(S(f_{i},g_{j})|G)=0$. Since $\overline{C(\Gamma_{1})}$ is arithmetically Cohen-Macaulay, $LM(f_{i})$ does not contain $x_{l}$, $NF(S(f_{i}, x^{b_{l}}_{l}-y^{a_{1}}_{1}\cdots y_{l}^{a_{l}}|G)=0$ and $ NF(g_{j}, x^{b_{1}}_{1}\cdots x^{b_{i}}_{i}-y^{a_{1}}_{1}|G)=0$ for  $ 1 \leq i\leq s$ and $ 1\leq j \leq t.$
\medskip

By Buchberger's criterion (\cite{BC},Theorem 1.7.3), this set is a Gr\"{o}bner  basis of $\frak{p}(\Gamma)$ with respect to the degree reverse lexicographic ordering $x_{1}>\dots>x_{l}>y_{1}>\dots>y_{k}$. By Lemma \ref{gbhom}, $ G^{h}=G_{1}^{h} \cup G_{2}^{h} \cup \{x^{b_{l}}_{l}-x_{0}^{b_{l}-(a_{1}+\cdots+a_{k})}y^{a_{1}}_{1}y^{a_{2}}_{2}\dots y^{a_{k}}_{k}\}$ is a Gr\"{o}bner basis of $\overline{\mathfrak{p}(\Gamma)}$ with respect to the degree reverse lexicographic ordering $x_{1} >\cdots>x_{l}>y_{1}>\cdots>y_{k}>x_{0}$. (Note that $G^{h}$ is a  homogenization of $G$ with respect to $x_{0}$.) \qed
\medskip
 
\begin{theorem}\label{GT}  Let $\Gamma_{1} = \Gamma(m_{1},\dots, m_{l})$ and $\Gamma_{2} = \Gamma(n_{1},\dots, n_{k})$ be two numerical semigroups with $m_{1} < \cdots < m_{l}$ and $n_{1} < \cdots < n_{k}$, take $p =b_{1}m_{1} +\cdots +b_{l}m_{l} \in \Gamma_{1}$ and $q = a_{1}n_{1} +\cdots +a_{k}n_{k} \in\Gamma_{2}$ be two positive integers satisfying
$\gcd(p, q) = 1 $ such that $\Gamma=\Gamma_{1}\star_{p,q} \Gamma_{2}$ is the star gluing of $\Gamma_{1}$ and $\Gamma_{2}$ . If the associated projective closures $\overline{C(\Gamma_{1})}$ and $\overline{C(\Gamma_{2})}$ are arithmetically Cohen-Macaulay, then the projective closure $\overline{C(\Gamma)}$ associated to $\Gamma$ is arithmetically Cohen-Macaulay.
\end{theorem}

\proof From Lemma \ref{GG}, $G^{h}=G_{1}^{h} \cup G_{2}^{h} \cup \{x^{b_{l}}_{l}-x_{0}^{b_{l}-(a_{1}+\cdots+a_{k})}y^{a_{1}}_{1}y^{a_{2}}_{2}\dots y^{a_{k}}_{k}\}$ is a Gr\"{o}bner bases of $\overline{\frak{p}(\Gamma)}$  with respect to the degree reverse lexicographic ordering $x_{1}>\dots>x_{l}>y_{1}>\dots>y_{k}>x_{0}$, where $G_{1}$ and $G_{2}$ are Gr\"{o}bner bases of $\frak{p}(\Gamma_{1})$ and $\frak{p}(\Gamma_{2})$ with respect to the degree reverse lexicographic ordering $x_{1}>\dots>x_{l}$ and $y_{1}>\dots>y_{k}$ respectively. Since $ \overline{C(\Gamma_{1})}$ and  $ \overline{C(\Gamma_{2})}$ are  arithmetically Cohen-Macaulay, by Lemma \ref{criteria},  $y_{k}$ does not divide leading monomial of any element in $G_{1}^{h}$ and $G_{2}^{h}$ and $LM(x^{b_{l}}_{l}-x_{0}^{b_{l}-(a_{1}+\cdots+a_{k})}y^{a_{1}}_{1}y^{a_{2}}_{2}\dots y^{a_{k}}_{k})=x_{l}^{b_{l}}$. Thus $y_{k}$ does not divide leading monomial of any element in $G^{h}$ which is a Gr\"{o}bner basis with respect to to the degree reverse lexicographic ordering $x_{1} >...>x_{l}>y_{1}>...>y_{k}>x_{0}$. Thus, by Theorem \ref{criteria}, $\overline{C(\Gamma)}$ is arithmetically Cohen-Macaulay. \qed
\medskip

\begin{theorem}\label{HIL}
A Noetherian graded domain $R$ over $k=R_{0}$ is Gorenstein if and 
only if it is Cohen-Macaulay and Hilbert series satisfy the identity 
$$H(R,1/t)=(-1)^{\mathrm{dim} R}\cdot t^{l}\cdot H(R,t)$$ for some $l \in \mathbb{Z}.$
\end{theorem}

\proof Follows from Theorem 4.4 in \cite{ST}. \qed
\medskip

\begin{lemma}\label{AP}
Let $\Gamma$ be a numerical semigroup and let $n$ be a positive integer of $\Gamma$. Let $Ap(\Gamma,n)=\{a_{0}<a_{1}<\cdots<a_{n-1}\}$ be the Ap\'{e}ry set of $n$ in $\Gamma$. Then $S$ is symmetric if and only if $a_{i}+a_{n-1-i}=a_{n-1}$ for all $i \in \{0,\dots,n-1\}$.
\end{lemma}

\proof See Proposition 4.10 in \cite{RS}. \qed
\medskip

\begin{lemma}\label{GC}
Assume that $\overline{C(\Gamma)}$ is Cohen-Macaulay. Let $ B=\{b_{i}=(\upsilon_{i},\mu_{i})\}$ be the Ap\'{e}ry set with respect to $n_{r}$ of 
$\overline{\Gamma}$, ordered so that $\upsilon_{1} < \upsilon_{2} <\dots<\upsilon_{n_{r-1}}$. Then $\overline{C(\Gamma)}$ is 
Gorenstein if and only if $b_{n_{r-1}}=b_{i}+b_{n_{r-1-i}}$, for all $i=1,\dots,n_{r-2}$.
\end{lemma}

\proof See Proposition 4.11 in \cite{CN}. \qed

\vspace{11pt}
\allowdisplaybreaks

\begin{theorem}\label{gorenstein} Let $\Gamma_{1} = \Gamma(m_{1},\dots, m_{l})$ and $\Gamma_{2} = \Gamma(n_{1},\dots, n_{k})$ be two numerical semigroups with $m_{1} < \cdots < m_{l}$ and $n_{1} < \cdots < n_{k}$. Let 
$p =b_{1}m_{1} +\cdots +b_{l}m_{l} \in \Gamma_{1}$ and $q = a_{1}n_{1} +\cdots +a_{k}n_{k} \in\Gamma_{2}$ be two positive integers satisfying
$\gcd(p, q) = 1 $ such that $\Gamma=\Gamma_{1}\star_{p,q} \Gamma_{2}$ is the star gluing of $\Gamma_{1}$ and $\Gamma_{2}$ . If the associated projective closures $\overline{C(\Gamma_{1})}$ and $\overline{C(\Gamma_{2})}$ are Gorenstein, then the projective closure $\overline{C(\Gamma)}$ associated to $\Gamma$ is Gorenstein.
\end{theorem}

\proof Since $\overline{C(\Gamma_{1})}$ and $\overline{C(\Gamma_{2})}$ are Gorenstein. From Lemma \ref{GD}, the corresponding Ap\'{e}ry set will be good, by Definition \ref{good} and Lemma \ref{GC}, Ap\'{e}ry set of $\Gamma_{1}$ and $\Gamma_{2}$ will satisfy the condition of Lemma \ref{AP}. They will be symmetric and gluing of two symmetric numerical semigroups is symmetric [prop 9.11, \cite{RS}], this implies $\Gamma$ is symmetric. Using Lemma \ref{DIG}, assume that $\frak{p}(\Gamma)=\langle \{f_{1},\dots,f_{s},g_{1},\dots,g_{t},x_{l}^{b_{l}}-y_{1}^{a_{1}}\dots y_{l}^{a_{l}}\}\rangle \subset k[x_{1},\dots,x_{l},y_{1},\dots,y_{k}](=:R)$ is the defining ideal of $C(\Gamma)$, where $\frak{p}(\Gamma_{1})=\langle\{f_{1},\dots,f_{s}\}\rangle \, \, \text{and}\,\,\frak{p}(\Gamma_{2})=\langle\{g_{1},\dots,g_{t}\}\rangle$. From Lemma \ref{GG}, $G= \{f_{1},\dots,f_{s},g_{1},\dots,g_{t},x_{l}^{b_{l}}-y_{1}^{a_{1}}\dots y_{l}^{a_{l}}\}$ is a Gr\"{o}bner basis of $\frak{p}(\Gamma)$ with respect to the degree reverse lexicographic ordering $x_{1}>\dots>x_{l}>y_{1}>\dots>y_{k}$ and $G^{h}$ is a Gr\"{o}bner basis of $\overline{\frak{p}(\Gamma)}\subset R[x_{0}]$ with respect to the degree reverse lexicographic ordering $x_{1}>\dots>x_{l}>y_{1}>\dots>y_{k}>x_{0}$. 
\medskip

From Theorem \ref{GT}, $\overline{C(\Gamma)}$ is arithmetically Cohen-Macaulay and note that $(in_{<_{rev}}(\frak{p}(\Gamma)))=(in_{<_{rev}}(\overline{\frak{p} (\Gamma)})$ with respect to the degree reverse lexicographic ordering $x_{1}>\dots>x_{l}>y_{1}>\dots>y_{k}>x_{0}$.. Hilbert function of ideal and its initial ideal is same, hence 
\begin{align*}
H(R[x_0]/\overline{\frak{p} (\Gamma)},t)=&H(R[x_0]/in_{<_{rev}}(\overline{\frak{p} (\Gamma)}),t)\\=&H(R/in_{<_{rev}}(\frak{p}(\Gamma)),t)\frac{1}{(1-t)}\\=&H(R/\frak{p}(\Gamma),t)\frac{1}{(1-t)}.
\end{align*}

Now we  apply Theorem \ref{HIL} and get,
\begin{align*}
H(R[x_0]/\overline{\frak{p} (\Gamma)},1/t)= 
& H(R[x_0]/in_{<_{rev}}(\overline{\frak{p} (\Gamma)}),1/t)\\
= & H(R/in_{<_{rev}}(\frak{p}(\Gamma)),1/t)\frac{1}{(1-\frac{1}{t})}\\
= & H(R/\frak{p}(\Gamma),1/t)\frac{1}{(1-\frac{1}{t})}\\
= & (-1)\cdot t^{l}\cdot H(R/\frak{p}(\Gamma),t)\frac{1}{(1-\frac{1}{t})} \,(\because R/\frak{p}(\Gamma)\,\, \text{is} \,\, \text{Gorenstein})\\
= & (-1)^{2}\cdot t^{l+1}\cdot H(R/\frak{p}(\Gamma),t)\frac{1}{(1-t)}\\=&(-1)^{2}\cdot t^{l+1}\cdot H(R[x_{0}]/\overline{\frak{p}(\Gamma)},t).
\end{align*}
Therefore, $\overline{C(\Gamma)}$  is Gorenstein. \qed
\medskip

\example\label{exm1}{\rm For $n \geq 2,$ let Arslan curve $\Gamma_{1}=\langle n(n+1),n(n+1)+1,(n+1)^{2},(n+1)^{2}+1\rangle$ and given integer $r \geq 3n+2 $, set $s=r(3n+2)+3 ,\, \text{Backelin curve}\,\, \Gamma_{2}= \langle s,s+3,s+3n+1,s+3n+2\rangle$. Note that $\overline{C(\Gamma_{1}})$ and $\overline{C(\Gamma_{2}})$ are arithmetically Cohen-Macaulay(see Theorem 2.8, \cite{TOA} and Proposition 5.2, \cite{CMC}).
\medskip

Choose $a\neq 1$ and $b\neq 1$ such that $p=b((n+1)^{2}+1), q= as, b >a,\,\, \text{and}\, (p,q)=1.$ Then $p$ and $q$ will satisfy the condition for the star gluing of $\Gamma_{1}$ and $\Gamma_{2}.$  Hence $\Gamma=\Gamma_{1}\star_{p,q} \Gamma_{2}=\langle q(n(n+1)),q(n(n+1)+1),q((n+1)^{2}),q(n+1)^{2}+1),ps,p(s+3),p(s+3n+1),p(s+3n+2)\rangle $. A Gr\"{o}bner basis of the defining ideal $\overline{\frak{p}(\Gamma)} \subset k[x_{0},x_{1},x_{2},x_{3},x_{4},y_{1},y_{2},y_{3},y_{4}]$ of $\overline{C(\Gamma)}$ with respect to the degree reverse lexicographic ordering $x_{1}>x_{2}>x_{3}>x_{4}>y_{1}>y_{2}>y_{3}>y_{4}>x_{0}$ is given by the binomials
 
$ f_{1}=\mathbf{x_{2}x_{3}^{3}}-x_{1}x_{4}^{3}$

 $ f_{(2,i)}=\mathbf{x_{1}^{n-i}x_{3}^{3i-1}}-x_{2}^{n-i+1}x_{4}^{3i-2}, 1 \leqslant i \leqslant n$
 
 $ f_{(3,j)}= \mathbf{x_{1}^{r-n+3+j}x_{2}^{n-1-j}}-x_{0}x_{3}^{2+3j}x_{4}^{r-1-3j}, 0 \leq j \leq n-1$

 $ f_{(4,j)}=\mathbf{x_{1}^{r-2n+3+j}x_{2}^{2n-j}}-x_{0}x_{3}^{3j+1}x_{4}^{r+1-3j}, 0 \leqslant j \leqslant n-1$
 
 $ f_{5}=\mathbf{x_{1}^{r-n+2}x_{2}^{n}x_{3}}-x_{0}x_{4}^{r+2}$
 
 $ f_{6}=\mathbf{x_{2}^{n+1}x_{3}}-x_{1}^{n}x_{4}^{2}$
 
 $ f_{7}=\mathbf{x_{2}^{2n+1}}-x_{1}^{2n-1}x_{3}x_{4}$
 
 $g_{i}=\mathbf{y_{1}^{n-i}y_{3}^{i+1}}-y_{2}^{n-i+1}y_{4}^{i},\text{for}\,\, 0 \leq i \leq n$
 
 $k_{i}=\mathbf{y_{1}^{i+1}y_{2}^{n-i}}-x_{0}y_{3}^{i}y_{4}^{n-i},\text{for} \,\, 0\leq i \leq n$
 
 $l= \mathbf{-y_{2}y_{3}}+y_{1}y_{4}$
 
 $l_{1}=\mathbf{x_{4}^{b}}-y_{1}^{a}$. 
  
Bold terms denote the leading monomial of corresponding binomials with respect to the degree reverse lexicographic ordering $x_{1}>x_{2}>x_{3}>x_{4}>y_{1}>y_{2}>y_{3}>y_{4}>x_{0}$.  
\medskip
 
$G_{1}^{h}=\{f_{1}^{h},f_{(2,i)}^{h},f_{(3,j}^{h},f_{4,j}^{h},f_{5}^{h},f_{6}^{h},f_{7}^{h} \vert 1 \leq i \leq n, 0 \leq j \leq n-1\}$ and $G_{2}^{h}= \{g_{i}^{h},k_{i}^{h},l \vert 0\leq i \leq n\}$ are Gr\"{o}bner bases of $\overline{\frak{p}(\Gamma_{1})} $ and $\overline{\frak{p}(\Gamma_{2})}$ with respect to the degree reverse lexicographic ordering $x_{1}>x_{2}>x_{3}>x_{4}>x_{0}$ and $y_{1}>y_{2}>y_{3}>y_{4}>x_{0}$ respectively(see \cite{BAC} and  \cite{CMC}). Let $G=G_{1}^{h} \cup G_{2}^{h} \cup \{l_{1}\}$. Since $\mathrm{gcd}(LM(f),LM(g),LM(l_{1}))=1,\,\text{for}\,f\in G_{1}^{h} \,\, \text{and}\, g \in G_{2}^{h}$, hence $spoly(f,g),spoly(f,l_{1})$ and $spoly(g,l_{1})$ reduce to $0$ when divided by $G$. Hence $G=G_{1}^{h} \cup G_{2}^{h} \cup l_{1}$ is a Gr\"{o}bner basis of the defining ideal $\frak{p}(\Gamma)$ with respect to the degree reverse lexicographic ordering $x_{1}>x_{2}>x_{3}>x_{4}>y_{1}>y_{2}>y_{3}>y_{4}>x_{0}$.  Since $y_{4}$ doesn't divide leading monomial of any element of Gr\"{o}bner basis $G$ of $\overline{\frak{p}(\Gamma)}$, from Theorem \ref{criteria}, $\overline{C(\Gamma)}$ is arithmetically Cohen-Macaulay.}
\medskip

\example{\rm Let $\Gamma_{1}=\Gamma(5,7,9,11)$ and $\Gamma_{2}=\Gamma(4,7,10)$ be two numerical semigroup. Since $\Gamma_{1}$ and $\Gamma_{2}$ are minimally generated by arithmetic progression, hence $\overline{C(\Gamma_{1})}$ and $\overline{C(\Gamma_{2})}$ are arithmetically Cohen-Macaulay (see Theorem 2.2, \cite{ACM}). 
\medskip

Consider the star gluing of $\Gamma_{1}$ and $\Gamma_{2}$ with  $p=2(11)=22$ and $q=4+7=17$, 
$\Gamma= \Gamma_{1}\star_{22,17}\Gamma_{2}=\langle 85,119,153,187,88,154,220\rangle.$ By the same argument as in Example \ref{exm1}, Gr\"{o}bner basis of the defining ideal $\overline{\frak{p}(\Gamma)} \subset k[x_{0},x_{1},x_{2},x_{3},x_{4},y_{1},y_{2},y_{3}]$ of $\overline{C(\Gamma)}$ with respect to the degree reverse lexicographic ordering $x_{1}>x_{2}>x_{3}>x_{4}>y_{1}>y_{2}>y_{3}>x_{0}$ is given by $f_{1}=\mathbf{x_{3}^2}-x_{2}x_{4},f_{2}=\mathbf{x_{2}x_{3}}-x_{1}x_{4},f_{3}=\mathbf{x_{2}^{2}}-x_{1}x_{3},f_{4}=\mathbf{x_{1}^{3}x_{3}}-x_{0}^{2}x_{4}^{2},f_{5}=\mathbf{x_{1}^{4}}-x_{0}^{2}x_{3}x_{4},f_{6}=\mathbf{y_{1}^{5}}-x_{0}^{3}y_{3}^{2},f_{7}=\mathbf{y_{2}^{2}}-y_{1}y_{3},f_{8}=\mathbf{x_{4}^{2}}-y_{2}y_{3}$. Since $y_{3}$ doesn't divide $\mathrm{LM}(f_{i})$ for all $1 \leq i \leq 8.$ From Theorem \ref{criteria}, $\overline{C(\Gamma)}$ is arithmetically Cohen-Macaulay.
\par Using Macaulay 2\cite{Macaulay}, the Hilbert series of $\overline{C(\Gamma)}$ is given by
\begin{align*}
(1-t)^{8}H(R[x_{0}]/\overline{\frak{p}(\Gamma)},t)=&t^{16}-5t^{14}+2t^{13}+5t^{12}-2t^{11}+t^{10}-4t^{8}+t^{6}\\&-2t^{5}+5t^{4}+2t^{3}-5t^{2}-1 \\=&(-1)^{2}\,t^{8}H(R[x_{0}]/\overline{\frak{p}(\Gamma)},\frac{1}{t})(1-t)^{8}
\end{align*}

$\Rightarrow H(R[x_{0}]/\overline{\frak{p}(\Gamma)},\frac{1}{t})=(-1)^{2}\,t^{-8}H(R[x_{0}]/\overline{\frak{p}(\Gamma)},t)$

By Theorem \ref{HIL}, $\overline{C(\Gamma)}$ is Gorenstein also.}
\medskip

\example{\rm: Consider $\Gamma_{1}=\Gamma(67,70,74,75)$ and $\Gamma_{2}=\Gamma(4,7,10)$. Note that $\overline{C(\Gamma_{1}})$ is arithmetically Cohen-Macaulay(see Theorem 2.8, \cite{BAC}) and since $\Gamma_{2}$ is minimally generated by arithemetic sequence, $\overline{C(\Gamma_{2}})$ is arithmetically Cohen-Macaulay.

Consider the star gluing of $\Gamma_{1}$ and $\Gamma_{2}$ with $p=2(75)=150$ and $q=4+7=11,\, \Gamma=\Gamma_{1} \star_{150,11} \Gamma_{2}=\langle 737,770,814,825,600,1050,1500 \rangle $. By the same argument as in Example \ref{exm1}, Gr\"{o}bner basis of the defining ideal $\overline{\frak{p}(\Gamma)} \subset k[x_{0},x_{1},x_{2},x_{3},x_{4},y_{1},y_{2},y_{3}]$ of $\overline{C(\Gamma)}$ with respect to the degree reverse lexicographic ordering $x_{1}>x_{2}>x_{3}>x_{4}>y_{1}>y_{2}>y_{3}>x_{0}$ is given by $f_{1}=\mathbf{x_{1}x_{3}^{2}}-x_{2}^{2}x_{4},f_{2}=\mathbf{x_{2}x_{3}^{3}}-x_{1}x_{4},f_{3}=\mathbf{x_{2}^{3}x_{3}}-x_{1}^{2}x_{4}^{2},f_{4}=\mathbf{x_{3}^{5}}-x_{2}x_{4}^{4},f_{5}=\mathbf{x_{2}^{5}}-x_{1}^{3}x_{3}x_{4},f_{6}=\mathbf{x_{1}^{9}x_{2}}-x_{3}^{2}x_{4}^{7}x_{0},f_{7}=\mathbf{x_{1}^{10}}-x_{2}x_{4}^{8}x_{0},f_{8}=\mathbf{x_{1}^{8}x_{2}^{2}x_{3}}-x_{4}^{10}x_{0},f_{9}=\mathbf{x_{1}^{7}x_{2}^{4}}-x_{3}x_{4}^{9}x_{0},f_{10}=\mathbf{x_{1}^{8}x_{2}^{3}}-x_{3}^{4}x_{4}^{6}x_{0},f_{11}=\mathbf{y_{1}^{5}}-x_{0}^{3}y_{3}^{2},f_{12}=\mathbf{y_{2}^{2}}-y_{1}y_{3},f_{13}=\mathbf{x_{4}^{2}}-y_{1}y_{2}$.
 Since $y_{3}$ doesn't divide $ LM(f_{i})$ for all $1 \leq i \leq 13.$ From Theorem \ref{criteria}, $\overline{C(\Gamma)}$ is arithmetically Cohen-Macaulay.}													
\bigskip

\section{Betti Sequence of Projective Closure}

Finding the minimal free resolution of modules and their invariants such as 
Betti numbers are the classical problem in commutative algebra. Usually it 
is very difficult to find the explicit description of the differential in 
a minimal free resolution, we can try to obtain information about the Betti 
numbers.  The $i^{th} \,Betti \, number$ of an $R-$ module $M$ is defined 
as 
$$\beta_{i}^{R}(M)=\mathrm{rank}(F_{i})$$ 
where $F_{i}$ is the free module appearing in the minimal free resolution of $M$, $F.$
$$\mathbf{F}: 0 \rightarrow F_{i} \rightarrow F_{i-1} \rightarrow \cdots \rightarrow F_{1} \rightarrow F_{0}.$$   

In the survey article on Betti numbers for affine monomial curves, written by D.Stamate \cite{STA}, 
one can find several useful information. However, the study of Betti numbers of the projective 
closure of affine monomial curves has not been done very extensively yet. 
In this paper our aim is to study the Betti numbers of the projective closure 
of  numerical semigroup ring and give a condition on a Gr\"{o}bner basis of 
the defining ideal of affine monomial curve, such that the Betti numbers of 
the projective closure are exactly the same as the Betti numbers of the 
affine monomial curve; see Theorem \ref{BSP}. Gimenez et al.( \cite{BAS}) have given 
the description of the Betti sequence of numerical semigroup ring 
generated by an arithmetic sequence. In Corollary \ref{BPAS}, we use 
Theorem \ref{BSP} to prove that for this class of affine monomial curves, 
the Betti numbers of the projective closure are the same as those of the affine curve. 
\medskip
 
\begin{lemma}\label{gbhom}
Let $I$ be an ideal in $R=k[x_{1},\ldots,x_{r}]$ and $I^{h}\subset R[x_{0}]$ its homogenization w.r.t. the variable $x_{0}$. Let $<$ be any reverse lexicographic monomial order on $R$ and $<_{0}$ the reverse lexicographic monomial order on $R[x_{0}]$ extended  from $R$ such that $x_{i}>x_{0}$ is the least variable.

If $\{f_{1},\ldots,f_{n}\}$ is the reduced Gr\"{o}bner basis for $I$ w.r.t $<$, then $\{f_{1}^{h},\ldots,f_{n}^{h}\}$ is the reduced Gr\"{o}bner basis for $I^{h}$ w.r.t $<_{0}$, and $\mathrm{in}_{<_{0}}(I^{h})=(\mathrm{in}_{<}(I))R[x_{0}]$.
\end{lemma}
\proof See Lemma 2.1 in \cite{CMC}.\qed
\medskip

\begin{theorem}\label{criteria}
Let $\mathbf{n} = (n_{1},\ldots,n_{r})$ be a sequence of positive integers with $n_{r}> n_{i}$ for all $i<n$. Let $<$ any reverse lexicographic order on $R=K[x_{1},\ldots,x_{r}]$ such that $x_{i}>x_{r}$ for all $1\leq i<r$. and $<_{0}$ the induced reverse lexicographic order on $R[x_{0}]$, where $x_{r}>x_{0}$. Then the following conditions are equivalent:
\begin{enumerate}
\item[(i)] The projective monomial curve $\overline{C(n_{1},\ldots,n_{r})}$ is arithmetically Cohen-Macaulay.
\item[(ii)] $\mathrm{in}_{<_{0}}((\mathfrak{p}(n_{1},\ldots,n_{r}))^{h})$ (homogenization w.r.t. $x_{0}$) is a Cohen-Macaulay ideal.
\item[(iii)] $\mathrm{in}_{<}(\mathfrak{p}(n_{1},\ldots,n_{r}))$ is a Cohen-Macaulay ideal.
\item[(iv)] $x_{r}$ does not divide any element of $G(\mathrm{in}_{<}(\mathfrak{p}(n_{1},\ldots,n_{r})))$.
\end{enumerate}
\end{theorem}
\proof See Theorem 2.2 in \cite{CMC}. \qed
\medskip

\begin{lemma} \cite{Peeva}\label{TOR}
Suppose that $U$ and $W$ are graded finitely generated $R-$ modules. Let $u$ be both $R-$ regular and $U-$ regular. Suppose that $uW=0.$
Then $$\mathrm{Tor}_{i}^{R}(U,W) \cong \mathrm{Tor}_{i}^{R/(u)}(U/uU,W)$$, for all $i \geq 0.$  
\end{lemma}														

\begin{lemma}[Corollary 10.61, \cite{Rotman}]\label{TOR BASE}
If $ _{R}B_{S}$ is flat on either side, then
$$ \mathrm{Tor}_{n}^{S}(A \otimes_{R} B,C) \cong \mathrm{Tor}_{n}^{R}(A,B \otimes_{S}C).$$
\end{lemma}	

\allowdisplaybreaks
\begin{remark}\label{BP}
Suppose $x_{n}$ is both $R$-regular and $\dfrac{R}{G}-$ regular, from Lemma \ref{TOR}, we have
\begin{align*}
\mathrm{Tor}_{n}^{R}\left(\dfrac{R}{G},k\right)\cong & \mathrm{Tor}_{n}^{R/x_{n}}\left(\frac{R}{(G,x_{n})},k\right)\\ 
\end{align*}
 from Lemma \ref{TOR BASE},
\begin{align*}
 \mathrm{Tor}_{n}^{R/(x_{n})}\left(\frac{R}{(G,x_{n})},k\right)=&\mathrm{Tor}_{n}^{R/(x_{n})}\left(\frac{R}{(G,x_{n})}\otimes_{R} \frac{R}{x_{n}},k\right)\\\cong & \mathrm{Tor}_{n}^{R}\left(\frac{R}{(G,x_{n})},\frac{R}{x_{n}} \otimes_{R/x_{n}}k\right) \\ \cong & \mathrm{Tor}_{n}^{R}\left(\frac{R}{(G,x_{n})},k\right).
\end{align*}		

Hence, Betti numbers are preserved under dividing by regular elements.
\end{remark}
\medskip

\begin{theorem}\label{BSP}
Let $\Gamma$ be a numerical semigroup such that $\overline{C(\Gamma)}$ is arithmetically Cohen-Macaulay. If there exist a Gr\"{o}bner basis $G$ with respect to the degree reverse lexicographic ordering $x_{i}>x_{0}$, of the defining ideal $\mathfrak{p}(\Gamma)$ of $C(\Gamma)$ such that $x_{n}$ belongs to the support of all nonhomogenous element of $G$, then $$\beta_{i}(C(\Gamma))=\beta_{i}(\overline{C(\Gamma)}).$$															
\end{theorem}																																	
\proof Let $G=\{f_{1},\dots,f_{t},g_{1},\dots,g_{r}\}$ be a Gr\"{o}bner basis of $\mathfrak{p}(\Gamma)$ with respect to the degree reverse lexicographic order $x_{i}>x_{0}$. From Theorem \ref{gbhom}, note that $G^{h}$ is a Gr\"{o}bner basis of $\overline{\mathfrak{p}(\Gamma)}$ with respect to the degree reverse lexicographic order $x_{i}>x_{0}$. Since $\overline{C(\Gamma)}$ is arithmetically Cohen-Macaulay then $x_{n}$ does not divide any leading monomial of $G$. Let $f_{1},\dots,f_{t}$ be homogeneous binomials and $g_{1},\dots,g_{r}$ be non-homogeneous elements. From our assumption, the terms of $g_{i}$ which are not the leading term is divisible by $x_{n}$, for $i=1,\dots,r$. 

Consider the natural map  
\[
\pi : S=k[x_{0},\dots,x_{n}] \longrightarrow k[x_{1},\dots,x_{n-1}].
\]

where $\pi(x_{n})=0,\pi(x_{0})=0$ and $\pi(x_{i})=x_{i}$ for $i=1,\dots,n-1$. Note that $$\frac{k[x_{0},\dots,x_{n}]}{(G^{h},x_{0},x_{n})} \cong \frac{k[x_{1},\dots,x_{n-1}]}{\pi(G^{h})}$$ as $S/x_{0}$-module. 

Since $\pi(g_{i})$ is a monomial which is not divisible by $x_{0}$, $x_{n}$ and any term of $f_{i}$ is not divisible by $x_{0}$. Hence $\pi(G^{h})=\pi(G).$

Consider
\begin{align*}
\beta_{i}^{S}\left(\frac{k[x_{0},\dots,x_{n}]}{G^{h}}\right)=&\beta_{i}^{S}\left(\frac{k[x_{0},\dots,x_{n}]}{(G^{h},x_{n})}\right)\\
\because \, x_{0} \, \text{is} \, S/(x_{n},G^{h})\text{-regular, by Lemma}\, \ref{TOR}\\ \beta_{i}^{S}\left(\frac{k[x_{0},\dots,x_{n}]}{(G^{h},x_{n})}\right) =& \beta_{i}^{S/x_{0}}\left(\frac{k[x_{0},\dots,x_{n-1}]}{(G^{h},x_{0},x_{n})}\right)\\=& \beta_{i}^{S/x_{0}}\left(\frac{k[x_{1},\dots,x_{n-1}]}{\pi(G^{h})}\right)\\=& \beta_{i}^{R}\left(\frac{k[x_{1},\dots,x_{n-1}]}{\pi(G)}\right)\\=&\beta_{i}^{R}\left(\frac{k[x_{1},\dots,x_{n}]}{(G,x_{n})}\right)\\ = &\beta_{i}^{R}\left(\frac{k[x_{1},\dots,x_{n}]}{G}\right)(\text{from Remark \ref{BP}) }.
\end{align*}

Hence $\beta_{i}(C(\Gamma))=\beta_{i}(\overline{C(\Gamma)})$ for all $ i \geq 1$. \qed						\medskip

\begin{corollary}\label{BPAS}
Let $\Gamma=\{m_{1},\dots,m_{n}\} $ be a numerical semigroup minimally generated by arithmetic sequence $m_{1}<m_{2}<\dots<m_{n}$ such that $m_{i}=m_{1}+(i-1)d, 1 \leq i \leq n$ and $m_{1}=q(n-1)+r, r \in [1,n]$, then $$\beta_{i}(C(\Gamma))=\beta_{i}(\overline{C(\Gamma)}).$$	
\end{corollary}		

\proof: Since $G=\{x_{i}x_{j}-x_{i-1}x_{j+1} \vert 2 \leq i \leq  j \leq n-1\} \cup \{x_{1}^{q+d}x_{i}-x_{r+i}x_{n}^{q} \vert 1 \leq i \leq n-r\}$ is a Gr\"{o}bner basis of the defining ideal $\mathfrak{p}(\Gamma)$ of $C(\Gamma)$ with respect to the degree reverse lexicographic ordering $x_{1}>x_{2}>\cdots>x_{n}$ and $ in_{<}(G)=\{x_{i}x_{j} \vert 2 \leq i \leq j \leq n-1\} \cup \{x_{1}^{q+d}x_{i} \vert 1 \leq i \leq n-r\}$( see \cite{Bermejo} ). Note that $x_{n}$ belongs to the support of all non-homogeneous elements of $G$. From Theorem \ref{BSP}, we have $\beta_{i}^{R}(C(\Gamma))=\beta_{i}^{S}(\overline{C(\Gamma)}).$					\qed																	\bigskip

\section{Betti sequence for simple gluing of numerical semigroup}
In this section, we are interested in the case when one of the glued semigroups 
is $\mathbb{N}$. We establish the relation between the Betti numbers of the affine 
monomial curve and the Betti numbers of the projective closure of the affine monomial 
curve defined by the numerical semigroup obtained by simple gluing of the numerical 
semigroup defing the affine monomial curve. First we recall the 
definition of simple gluing.

\begin{definition}[\cite{SG}]
Let $\Gamma$ be a numerical semigroup
minimally generated by $a_{1} < \dots < a_{n},$ and $c > 1$ and $d$ coprime integers such that $d \in  \Gamma \setminus \{a_{1}, \dots , a_{n}\}$. Then
$$\Gamma_{1} = \langle c\Gamma,d \rangle = \langle ca_{1}, \dots , ca_{n}, d \rangle $$ is said to be obtained from $\Gamma$ by a simple gluing.

\end{definition}

The order of $d$ in $\Gamma$ is defined as $$ord_{\Gamma}(d)=max\{\sum_{i=1}^{n}\lambda_{i}:d=\sum_{i=1}^{n}\lambda_{i}a_{i},\lambda_{i} \in \mathbb{N}\}.$$

\begin{lemma}(\cite{Watanabe}, Lemma 1)\label{DSG} Let $\Gamma$ be a numerical semigroup minimally generated by $a_{1} < \dots < a_{n},$ and $c > 1$ and $d$ coprime integers such that $d \in  \Gamma \setminus \{a_{1}, \dots , a_{n}\}$. Suppose $\Gamma_{1} = \langle c\Gamma,d \rangle = \langle ca_{1}, \dots , ca_{n}, d \rangle $ is simple gluing of $\Gamma$. If we write $d = \lambda_{1}a_{1} + \dots + \lambda_{n}a_{n}$ with $\lambda_{i}$ non-negative integers and $\sum_{i=1}^{n}\lambda_{i}$
maximal, then defining ideal $\frak{p}(\Gamma_{1})$ is generated by $\frak{p}(\Gamma)$ and $f=x_{n+1}^{c}-\prod_{i=1}^{n}x_{i}^{\lambda^{i}}$.
\end{lemma}

\begin{lemma}\label{gbs}

Let $\Gamma$ be a numerical semigroup minimality generated by $a_{1}<\dots<a_{n}$. Suppose $c>1$ and $d$ are co-prime integers with $d\in \Gamma \setminus\{a_{1},\dots,a_{n}\}$. Let $\Gamma_{1}= \langle c\Gamma,d \rangle $ is a simple gluing.

\begin{enumerate}[(a)]

\item Let $c \geq ord_{\Gamma}(d)$ and $G$ be a Gr\"{o}bner basis of $\frak{p}_{\Gamma}$ with respect to the degree reverse lexicographic monomial ordering $<$ such that $x_{n+1}$ is the largest variable then $(G,f)$ will be a Gr\"{o}bner basis of $\frak{p}_{\Gamma_{1}}$ with respect to the degree reverse lexicographic ordering $<$, where $f=x_{n+1}^{c}-\prod_{i=1}^{n}x_{i}^{\lambda^{i}}$ with  $d = \lambda_{1}a_{1} + \dots + \lambda_{n}a_{n}$ with $\lambda_{i}$ non-negative integers and $\sum_{i=1}^{n}\lambda_{i}$
maximal.
\item let $c \leq ord_{\Gamma}(d)$ and $G$ be a Gr\"{o}bner basis of $\frak{p}_{\Gamma}$ with respect to the negative degree reverse lexicographic monomial ordering $<$ such that $x_{n+1}$ is the largest variable then $(G,f)$ will be the Gr\"{o}bner basis of $\frak{p}_{\Gamma_{1}}$ with respect to the negative degree reverse lexicographic monomial ordering $<$, where $f=x_{n+1}^{c}-\prod_{i=1}^{n}x_{i}^{\lambda^{i}}$ with  $d = \lambda_{1}a_{1} + \dots + \lambda_{n}a_{n}$ with $\lambda_{i}$ non-negative integers and $\sum_{i=1}^{n}\lambda_{i}$ maximal.
\end{enumerate}
\end{lemma}

\proof(a): Let $G=\{g_{1},\dots,g_{r}\}$ be a Gr\"{o}bner basis with respect to the degree reverse lexicographic monomial ordering $<$ such that $x_{n+1}$ is the largest and $LM_{<}(f)=x_{n+1}^c$. Since leading term of any $g_{i}$ does not contain $x_{n+1}$,  $\mathrm{gcd}(LM_{<}(f),LM_{<}(g_{i}))=1$. For any $i, spoly(f,g_{i})\rightarrow_{G} 0$, hence $(G,f)$ will be a Gr\"{o}bner basis of $\frak{p}_{\Gamma_{1}}$ with respect to the degree reverse lexicographic monomial ordering $<$.

\proof(b): Since leading monomial of $f$ with respect to the degree reverse lexicographic monomial ordering $<$ such that $x_{n+1}$ is the largest, is again $x_{n+1}^{c}$.  Therefore $(G,f)$ will be the Gr\"{o}bner basis of $\frak{p}_{\Gamma_{1}}$ with respect to the negative degree reverse lexicographic monomial ordering $<$.\qed

\begin{definition}\label{•}
If $\mu:F \rightarrow G$ is a map of complexes, and we write $\phi$ and $\psi$, respectively, for the differentials of $F$ and $G$, then the mapping cone $M(\mu)$ of $\alpha$ is the complex such that $M(\mu)_{i}=F_{i-1} \oplus G_{i}$, with differential $\delta$ defined as follows 
\[
\mu_{\vert_{F_{i-1}}}=-\phi_{i-1}+\mu_{i}=F_{i-1} \rightarrow F_{i-2} \oplus G_{i-1}
\]
\[
\mu_{\vert_{G_{i}}}=\psi_{i}=G_{i} \rightarrow G_{i-1}
\]

\end{definition}

\begin{remark}\label{mp}Let $I$ be an ideal in $R$ and take an element $z\in R$. Then,$$0\rightarrow R/(I:z)\xrightarrow\mu R/I \rightarrow R/I+(z)\rightarrow 0 $$ is exact, where 
$\mu$ is the map given by multiplication by $z$. Now if $F$ resolves $R/(I:z)$, $G$ resolves $R/I$,and $\mu:F \rightarrow G$ is a map of complexes induced by $\mu$, then $M(\mu)$ resolves $R/I+(z)$.
\end{remark}

\medskip

\begin{theorem}\label{TP}
Let $\Gamma$ be a numerical semigroup minimally generated by $a_{1}<\dots < a_{n}$. Suppose $\Gamma_{1}=\langle c\Gamma,d \rangle$ be a numerical semigroup obtained from $\Gamma$ by simple gluing, where $c>1$ and $d$ are co-prime integers with $d\in \Gamma \setminus\{a_{1},\dots,a_{n}\}$. Then the following holds:
\begin{itemize}
\item[(i)] $\beta_{i}(C(\Gamma_{1}))=\beta_{i}(C(\Gamma))+\beta_{i-1}(C(\Gamma))$
\item[(ii)] $\beta_{i}(\overline{C(\Gamma_{1}}))=\beta_{i}(\overline{C(\Gamma}))+\beta_{i-1}(\overline{C(\Gamma}))$.
\end{itemize}
\end{theorem}

\proof 

\begin{enumerate}[(i)]

\item See Theorem 5.2(i) in \cite{STA}.

\item 
Let $F_{\Gamma}$ be a minimal free resolution of $\frac{R}{\overline{\frak{p}_{\Gamma}}}$. Since  $\overline{\frak{p}_{\Gamma}}$ is prime ideal, this implies $(\overline{\frak{p}_{\Gamma}}:f^{h})=\overline{\frak{p}_{\Gamma}}$,  $F_{\Gamma}$ will also resolve $R/(\overline{\frak{p}_{\Gamma}}:f^{h})$ minimally.
Now consider the complex map $\mu$ on $F_{\Gamma} \rightarrow F_{\Gamma}$ induced from multiplication by $f^{h}$. All the maps in the complex are also multiplication by $f^{h}$ and hence this is a complex map. From Remark \ref{mp} $M(\mu)$ resolve $R/(\overline{\frak{p}_{\Gamma}}+f^{h})=R/(\overline{\frak{p}_{\Gamma}},f^{h})$. Since $(\overline{\frak{p}_{\Gamma}},f^{h})$ is generating set of ideal $\overline{\frak{p}_{\Gamma_{1}}}$, $M(\mu)$ resolve $R/\overline{\frak{p}_{\Gamma_{1}}}$. Finally, since the maps in the complex map $\mu$ on $F_{\Gamma} \rightarrow F_{\Gamma}$ are  all multiplication by $f^{h}$ and $f^{h} \in (x_{0},x_{1},\dots.x_{n},x_{n+1})$. We see that $\delta.$ are the differentials in the $M(\mu)$, and for every $i$, one has $$\delta_{i}(M(\mu)_{i}) \subset (x_{0},x_{1},\dots x_{n},x_{n+1})M(\mu)_{i-1}$$ and  $M(\mu)$ resolves $R/\overline{\frak{p}_{\Gamma_{1}}}$ minimally and it give minimal free resolution of $\overline{C(\Gamma_{1}})$. Thus, by minimal free resolution $M(\mu)$ of $\overline{C(\Gamma_{1}})$, we have $\beta_{i}(\overline{C(\Gamma_{1}})=\beta_{i}(\overline{C(\Gamma}))+\beta_{i-1}(\overline{C(\Gamma}))$.
Moreover $\mathrm{pd}_{R}(\overline{C(\Gamma_{1}}))= \mathrm{pd}_{R}(\overline{C(\Gamma}))+1$.  \qed

\end{enumerate}

\begin{theorem}\label{aub}\cite{Peeva}
Let $V$ be a graded finitely generated $S$-module. Its projective dimension is $$\mathrm{pd}_{S}(V)=n-\mathrm{depth}(V).$$
\end{theorem}

\begin{theorem}
Let $\Gamma$ be a numerical semigroup minimality generated by $a_{1}<\dots<a_{n}$.  Suppose $\Gamma_{1}=\langle c\Gamma,d \rangle $ is a simple gluing, where $c>1$ and $d$ are co-prime integers with $d\in \Gamma\setminus\{a_{1},\dots,a_{n}\}$. If projective closure of affine curve $\overline{C(\Gamma})$ is arithmetically Cohen-Macaulay then $\overline{C(\Gamma_{1}})$ is arithmetically Cohen-Macaulay. 
\end{theorem}

\proof: Let $e$ be a projective dimension of $\overline{C(\Gamma})$. Since $\overline{C(\Gamma})$
is Cohen-Macaulay, $\mathrm{depth}(\overline{C(\Gamma}))=\mathrm{dim}(\overline{C(\Gamma}))=2$.
Now by Remark \ref{mp}, $\mathrm{pd}(\overline{C(\Gamma_{1}}))=e+1$. We will use \ref{aub}(Auslander-Buschsbaum), 
\begin{align*}
\mathrm{depth}(\overline{C(\Gamma_{1}}))=&(n+2)-\mathrm{pd}(\overline{C(\Gamma_{1}}))\\=&(n+2)-(e+1)\\=&(n+1)-e\\=&\mathrm{depth}(\overline{C(\Gamma}))\\=&2=\mathrm{dim}(\overline{C(\Gamma_{1})})
\end{align*}

Therefore, $\overline{C(\Gamma_{1})}$ is arithmetically Cohen-Macaulay.\qed

\smallskip 

\begin{corollary}
Let $\Gamma$ be a numerical semigroup minimality generated by $a_{1}<\dots<a_{n}$. Suppose $\Gamma_{1}=\langle c\Gamma,d \rangle $ is a simple gluing, where $c>1$ and $d$ are co-prime integers with $d\in \Gamma\setminus\{a_{1},\dots,a_{n}\}$ then $\overline{C(\Gamma_{1}})$ and $\overline{C(\Gamma})$ have same type.

In particular, if $\overline{C(\Gamma})$ is Gorenstein then $\overline{C(\Gamma_{1}})$ is Gorenstein. 
\end{corollary}																

\proof: Note that $n=\text{emb dim}(\Gamma)$.
From Theorem \ref{TP}, the type of $\overline{C(\Gamma_{1}})$ is $\beta_{n}(C(\Gamma_{1}))=\beta_{n}(C(\Gamma))+\beta_{n-1}(C(\Gamma))= \beta_{n-1}(C(\Gamma))$. \\
If  $\overline{C(\Gamma})$ is Gorenstein then 	
 $\overline{C(\Gamma})$ is arithmetically Cohen-Macaulay and type is 1. Now type of  $\overline{C(\Gamma_{1}})$ is $\beta_{n}(C(\Gamma_{1}))=\beta_{n}(C(\Gamma))+\beta_{n-1}(C(\Gamma))= \beta_{n-1}(C(\Gamma))=1$, and from Theorem \ref{TP}, $\overline{C(\Gamma_{1}})$ is arithmetically Cohen-Macaulay. Hence $\overline{C(\Gamma_{1}})$ is Gorenstein. \qed
\bigskip

\section{The type of projective closures of affine monomial curves}								
In this section, we will construct a numerical semigroup independent of the  embedding dimension such that the corresponding affine curve and its projective closure both are Goorenstein, that means they are arithmetically Cohen-Macaulay and their Cohen-Macaulay type is 1. 
\par Similarly, for a fixed positive integer $n>1$, we have constructed a family of numerical semigroups of arbitrary embedding dimension using simple gluing, such that the associated affine curve and its projective closure have the same Cohen-Macaulay type which is $n$.						

\begin{remark}[Remark 2.4 in \cite{Bermejo}]\label{CA1}
Let $m_{1}<m_{2}<\dots<m_{n}$ are in arithmetic progression such that $m_{1}=q(n-1)+r, \text{for}\, q \in \mathbb{N}, r \in \{1,\dots,n-1\}$. Consider $\Gamma= \langle m_{1},\dots,m_{n} \rangle$, then the first Betti number of minimal graded free resolution of $C(\Gamma)$ and 
$\overline{C(\Gamma})$ is $ {(n-1) \choose 2} + K$, where $K=n-r$.
\end{remark}

\smallskip

\begin{theorem}
Let $m_{1}<m_{2}<\dots<m_{n}$ are in arithmetic progression such that $m_{1} \equiv 2(n-1)$ and $\Gamma_{n}=\langle m_{1},\dots,m_{n} \rangle$. For $e>n$ construct a curve $\Gamma_{e}=\langle c\Gamma_{e-1},d \rangle$ such that $(c,d)=1$ and $d\in \Gamma_{n}$,  then $C(\Gamma_{e})$ and $\overline{C(\Gamma_{e}})$ have type  $1$ and their first Betti number is $\frac{n(n-1)-2}{2}+(n-e)$.

\end{theorem}

\proof: In \cite{Bermejo}, it is proved that  $\overline{C(\Gamma_{n}})$ is arithmetically Cohen-Macaulay and its Cohen Macaulay type is $1$.  By Theorem \ref{TP}, type of $\overline{C(\Gamma_{n+1}})$  is $1$. Inductively $\Gamma_{e}$ is a simple gluing of $\Gamma_{e-1}$,  again using Theorem \ref{TP} type of $\overline{C(\Gamma_{e}})$  is $1$. In \cite{TA}, it is proved that $C(\Gamma_{n})$ is Gorenstein. Thus $C(\Gamma_{e})$ and $\overline{C(\Gamma_{e}})$ have type $1$.

\par By Remark \ref{CA1}, the first Betti number of $\Gamma_{n}$ is $\frac{n(n-1)-2}{2}$. Using induction and from Theorem \ref{TP}, $\beta_{1}(\overline{C(\Gamma_{e}}))= \beta_{1}(\overline{C(\Gamma_{e-1}}))+1=\frac{n(n-1)-2}{2}+(n-e)$, which is same for $\beta_{1}(C(\Gamma_{e}))$.\qed
\medskip

\begin{theorem}
Let $n$ be any odd positive  integer greater than $2$, put $t=\frac{n+1}{2}$, $\Gamma_{4}=\langle t(t+1),t(t+1)+1,(t+1)^2,(t+1)^2+1 \rangle $. For $e>4$ construct a curve $\Gamma_{e}=\langle c\Gamma_{e-1},d \rangle$ such that $(c,d)=1$ and $d\in \Gamma$,  then $C(\Gamma_{e})$ and $\overline{C(\Gamma_{e})}$ have same type $n$. Their first Betti number is $ n+e-1 $ and $n+e$ respectively.
\end{theorem}											

\proof: In \cite{TOA}, it is proved that  $\overline{C(\Gamma_{4}})$ is arithmetically Cohen-Macaulay and its Cohen Macaulay type is $2t-1=n$.  By Theorem \ref{TP}, type of $\overline{C(\Gamma_{5}})$  is $n$. Inductively $\Gamma_{e}$ is a simple gluing of $\Gamma_{e-1}$ and using Theorem \ref{TP} type of $\overline{C(\Gamma_{e}})$  is $n$. The
semigroup ring $C(\Gamma_4)$ has the Betti sequence $(1,2t+2,4t,2t-1)$ (See \cite{STA}). Hence, $C(\Gamma_{e})$ and $\overline{C(\Gamma_{e}})  $ have same type $n$.
\par In \cite{TOA}, it is proved that the first Betti number of $\overline{C(\Gamma_{4})}$ is $ n+4$. Using induction and from Theorem \ref{TP}, $\beta_{1}(\overline{C(\Gamma_{e}}))= \beta_{1}(\overline{C(\Gamma_{e-1}}))+1=n+4+(e-5)+1=n+e$. Similarly, in \cite{STA} it is proved that the first Betti number of $C(\Gamma_{4})$ is $ n+3$. By same fact $\beta_{1}(C(\Gamma_{e}))=n+e-1$.\qed

\medskip
\begin{theorem}\textbf{(Gastinger)}\quad \label{gastinger} 
Let $A = k[x_{1},\ldots,x_{r}]$ be the polynomial ring, $I\subset A$ 
the defining ideal of a monomial curve defined by natural numbers 
$a_{1},\ldots,a_{r}$, whose greatest common divisor is $1$.  
Let $J$ be an ideal contained in $I$. Then $J = I$ if and 
only if $\mathrm{dim}_{k} A/\langle J + (x_{i}) \rangle =a_{i}$, 
for some $i$; equivalently 
$\mathrm{dim}_{k} A/\langle J + (x_{i}) \rangle =a_{i}$ for any $i$.
\end{theorem}  

\proof See in \cite{g}.\qed
\medskip

Let $I\subset k[x_{1},\ldots,x_{r}] $ be a monomial ideal, then it has unique minimal generating set. We denote the minimal generating set by $G(I)$.

\begin{theorem}\label{DC1} The defining ideal $\mathfrak{p}(\Gamma)$ of monomial
curve associated to $\Gamma=\langle t^2-t,t^2-t+1,t^2-1,t^2 \rangle $ is minimally generated by following binomials
\begin{itemize}
\item[(i)] $f_{(1,i)}=x_{3}^{2+i}x_{4}^{t-2-i}-x_{1}^{3+i}x_{2}^{t-2-i},\text{for}\,\, 0 \leq i \leq t-4$
\item[(ii)] $ f_{2}=x_{2}x_{3}-x_{1}x_{4}$
\item[(iii)] $f_{3}= x_{4}^{t-1}-x_{1}^{t}$
\item[(iv)] $f_{4}=x_{3}^{t-1}-x_{2}x_{4}^{t-2}$
\item[(v)] $f_{5}=x_{2}^{t-1}-x_{1}^{t-2}x_{3}$
\item[(vi)] $f_{(6,i)}=x_{1}^{i+1}x_{3}^{t-2-i}-x_{2}^{2+i}x_{4}^{t-3-i},\text{for}\,\, 0 \leq i \leq t-4$
\end{itemize}

\end{theorem}

\proof We proceed by Gastinger's Theorem stated as in Theorem \ref{gastinger}. Let $J=<\{f_{(1,i)},f_{2},f_{3},f_{4},f_{5},f_{(6,i)}\vert 0 \leq i \leq t-4 \}>$ and $J_{1}=\{x_{1}^{3+i}x_{2}^{t-2-i},x_{2}x_{3}$\\ $,x_{1}^{t},x_{3}^{t-1},x_{2}^{t-1}-x_{1}^{t-2}x_{3},x_{1}^{i+1}x_{3}^{t-2-i} \vert 0 \leq i \leq t-4 \}$. Then $J+(x_{4})=\langle J_{1} \rangle$. At first we note that $ J \subset \mathfrak{p}(\Gamma)$.

\medskip
$A/J+\langle x_{4}\rangle$ is the vector space over $k$ whose basis consists of the images
of monomials are listed below

\begin{itemize}
\item  $ T_{1}=\{x_{1}^{\alpha} : 0 \leq \alpha \leq (t-1)\}$
\item $ T_{2}=\{x_{1}^{\alpha}x_{2}^{\beta}: 1 \leq \alpha \leq 2, 1 \leq \beta \leq (t-2)\}$
\item $S_{2i}=\{x_{1}^{i+2}x_{2}^{\beta}: 1 \leq \beta \leq t-2-i\}, for \,\, 1 \leq i \leq t-3$
\item $S_{3i}=\{x_{1}^{i}x_{2}^{\beta}:1 \leq \beta \leq t-2-i\}, for \,\, 1 \leq i \leq t-3$
\item $ T_{3}=\{x_{2}^{\beta}:1 \leq \beta \leq t-2\}$

\item $ T_{4}= \{x_{3}^{\gamma}:1 \leq \gamma \leq t-2\}$
\item $ T_{5}=\{x_{1}^{t-2}x_{3},x_{1}^{t-1}x_{3}\}$
\end{itemize}
Then cardinality of this basis is,
\begin{align*}
\displaystyle \sum_{i=1}^{5}\vert T_{i}\vert + \displaystyle \sum_{i=1}^{t-3}\vert S_{2i} \vert +\displaystyle \sum_{i=1}^{t-3}\vert S_{3i} \vert =&t+2(t-2)+(t-2)+(t-2)+ 2\\ +&\displaystyle \sum_{i=0}^{t-4}(t-3-i)+\displaystyle \sum_{i=0}^{t-4}(t-3-i)\\ =&t^{2}
\end{align*}
Hence $dim_{k} \left(A/J+\langle x_{4}\rangle\right) =t^{2}$. By Gastinger Theorem, we have $ J=\mathfrak{p}(\Gamma)$. \qed

\medskip

\begin{theorem}\label{gbbc}
Let us consider the degree reverse lexicographic monomial order on $k[x_{1},\ldots,x_{4}] $ induced by $ x_{1}>x_{2}>x_{3}> x_{4}$ then $G=\{f_{(1,i)},f_{2},f_{3}$,\\$f_{4},f_{5},f_{(6,j)}\vert 0 \leq i,j \leq t-4 \}$ is Gr\"{o}bner basis of defining ideal $\mathfrak{p}(\Gamma)$ with respect to this order.
\end{theorem}
\proof We proceed by Buchberger's algorithm, 

\begin{enumerate}

\item
\begin{align*}
S(f_{(1,i},f_{(1,j)})=&x_{1}^{3+j}x_{2}^{t-2-j}x_{4}^{j-i}-x_{1}^{3+i}x_{3}^{j-i}x_{2}^{t-2-i}\\=&-(x_{2}x_{3}-x_{1}x_{4})(\displaystyle \sum_{l=0}^{j-i-1}x_{1}^{3+i+l}x_{2}^{t-3-i-l}x_{3}^{j-i-1-l}x_{4}^{l}).
\end{align*}

\item 

\begin{align*} 
S(f_{(1,i)},f_{2})=&x_{1}x_{3}^{1+i}x_{4}^{t-1-i}-x_{1}^{3+i}x_{2}^{t-1-i}, 0 \leq i \leq t-4 \\ =& x_{1}x_{3}f_{3}-x_{1}^{3}f_{5}, \quad \mathrm{for}\,\, i=0\\=& x_{1}f_{(1,i-1)}.\quad \mathrm{for}\,\, 1 \leq i \leq t-4\\
\end{align*} 

\item 
\begin{align*}
S(f_{(1,i)},f_{3})=&x_{3}^{2+i}x_{1}^{t}-x_{1}^{3+i}x_{2}^{t-2-i}x_{4}^{1+i},0 \leq i \leq t-4 
\\=& f_{(6,t-4)}x_{1}^{3}x_{3}^{i}+f_{2} \displaystyle \sum_{l=0}^{i-1}x_{1}^{3+l}x_{2}^{t-3-l}x_{3}^{i-1-l}x_{4}^{l+1}.
\end{align*}

\item
\begin{align*}
S(f(1,i),f_{4})=&x_{2}x_{4}^{2t-4-i}-x_{1}^{3+i}x_{2}^{t-2-i}x_{3}^{t-3-i}, 0 \leq i \leq t-4\\=&x_{2}x_{4}^{t-3-i}f_{3}-f_{2}(\displaystyle \sum_{l=0}^{t-4-i}x_{1}^{3+l+i}x_{2}^{t-3-l-i}x_{3}^{t-4-l-i}x_{4}^{l}). 
\end{align*}

\item

$\text{Since} \gcd(Lm(f_{(1,i)}),Lm(f_{5}))=1, \text{hence}\,  S(f_{(1,i)},f_{5}) \longrightarrow_{G} 0. $

\item 
We will rewrite $f_{(1,i)}=x_{3}^{t-i-2}x_{4}^{i+2}-x_{1}^{t-i-1}x_{2}^{i+2}$, for $0 \leq i \leq t-4$  \\
$ f_{(6,j)}=x_{1}^{l+j}x_{3}^{t-j-2}-x_{2}^{2+j}x_{4}^{t-3-j},$ for $0 \leq i \leq t-4$.

\begin{align*}
S(f_{(1,i),(6,j)}=&x_{2}^{2+j}x_{4}^{t-1-j+i}+-x_{1}^{t-i+j}x_{2}^{i+2}x_{3}^{i-j}\\=&f_{3}(x_{4}^{l-j}x_{2}^{2+j})-f_{2}\displaystyle\sum_{l=0}^{i-j-1}x_{1}^{t-i+j+l}x_{2}^{1+i-l}x_{3}^{i-j-1-l}x_{4}^{l},\text{for}\, i>j.
\\S(f_{(1,i),(6,j)}=&x_{2}^{2+j}x_{3}^{j-i}x_{4}^{t-1-j+i}-x_{1}^{j+t-i}x_{2}^{i+2}\\=&f_{6}(x_{1}^{j-i}x_{2}^{i+2})+f_{2}\displaystyle \sum_{l=0}^{j-i-1}x_{1}^{l}x_{2}^{1+j-l}x_{3}^{j-i-1-l}x_{4}^{t-j+i+l-1},\text{for}\, i<j.
\\S(f_{(1,i),(6,i)}=&x_{2}^{2+i}x_{4}^{t-1}-x_{1}^{t}x_{2}^{i+2}\\=&f_{6}x_{4}^{j+2}.
\end{align*}

\item $\text{Since} \gcd(Lm(f_{(2)}),Lm(f_{3}))=1, \text{hence}\,  S(f_{(2)},f_{3}) \longrightarrow_{G} 0. $

\item \begin{align*}
S(f_{2},f_{4})=& x_{2}^{2}x_{4}^{t-1}-x_{1}x_{3}^{t-2}x_{4}\\=&-x_{4}f_{(6,0)}. 
\end{align*}

\item \begin{align*}
S(f_{2},f_{5}=&x_{1}^{t-2}x_{3}^{2}-x_{1}x_{2}^{t-2}x_{4}\\=&x_{1}f_{(6,t-4)}.
\end{align*}

\item \begin{align*}
S(f_{2},f_{(6,i)}=&x_{2}^{3+i}x_{4}^{t-3-i}x_{4}-x_{1}^{i+2}x_{3}^{t-3-i}x_{4},0 \leq i \leq t-4\\ =& -x_{4}f_{(6,i+1)}, \quad \mathrm{for}\,\, 1 \leq i \leq t-3\\ =&
-x_{4}f_{5}\quad \mathrm{for}\,\, i=t-4.\\
\end{align*}

\item $\text{Since} \gcd(Lm(f_{3}),Lm(f_{4}))=1, \text{hence}\,  S(f_{3},f_{4}) \longrightarrow_{G} 0. $

\item $\text{Since} \gcd(Lm(f_{3}),Lm(f_{5}))=1, \text{hence}\,  S(f_{3},f_{5}) \longrightarrow_{G} 0. $

\item $\text{Since} \gcd(Lm(f_{3}),Lm(f_{(6,i)}))=1, \text{hence}\,  S(f_{3},f_{(6,i)}) \longrightarrow_{G} 0. $

\item $\text{Since} \gcd(Lm(f_{4}),Lm(f_{5}))=1, \text{hence}\,  S(f_{4},f_{5}) \longrightarrow_{G} 0. $

\item \begin{align*}
S(f_{4},f_{(6,i)}=&x_{2}^{2+i}x_{3}^{i+1}x_{4}^{t-3-i}-x_{1}^{i+1}x_{2}x_{4}^{t-2} \\=& f_{2}(\displaystyle \sum_{l=0}^{i}x_{1}^{l}x_{2}^{1+i-l}x_{3}^{i-l}x_{4}^{t-3-i+l}).
\end{align*}

\item $\text{Since} \gcd(Lm(f_{5}),Lm(f_{(6,i)}))=1, \text{hence}\,  S(f_{5},f_{(6,i)}) \longrightarrow_{G} 0. $

\item \begin{align*}
S(f_{(6,i)},f_{(6,j)}=&x_{3}^{j-i}x_{2}^{2+j}x_{4}^{t-3-j}-x_{1}^{j-i}x_{2}^{2+i}x_{4}^{t-3-i}\\=& f_{2}(\displaystyle \sum_{l=0}^{j-i-1}x_{1}^{l}x_{2}^{1+j-l}x_{3}^{j-i-1-l}x_{4}^{t-3-j+l} )
\end{align*}

\end{enumerate}
By Buchberger's Criterion  S-polynomials of all  generator reduces to zero, hence $G_{nr}$ is Gr\"{o}bner basis of $\mathfrak{p}(\Gamma)$ with respect to the degree reverse lexicographic monomial order $>$ induced by $ x_{1}>x_{2}>x_{3}> x_{4}$.\qed

\begin{theorem}
Let $n$ be any even positive  integer greater than $2$, put $t=\frac{n+4}{2}$, $\Gamma_{4}=\langle t^2-t,t^2-t+1,t^2-1,t^2 \rangle $. For $e>4$ construct a curve $\Gamma_{e}=\langle c\Gamma_{e-1},d \rangle$ such that $(c,d)=1$ and $d\in \Gamma$,  then $C(\Gamma_{e})$ and $\overline{C(\Gamma_{e})}$ have same type $n$ and same first Betti number $n-2+e$.
\end{theorem}
\proof: In \cite{CN}, it is proved that $\overline{C(\Gamma_{4}})$ is arithmetically Cohen-Macaulay and its Cohen Macaulay type is $2t-4=n$.  By Theorem \ref{TP}, type of $\overline{C(\Gamma_{5}})$  is $n$. Inductively $\Gamma_{e}$ is a simple gluing of $\Gamma_{e-1}$, using Theorem \ref{TP} type of $\overline{C(\Gamma_{e}})$  is $n$. The type of 
semigroup ring $C(\Gamma_4)$ is $2t-4=n$ (See \cite{CN}). Hence $C(\Gamma_{e})$ and $\overline{C(\Gamma_{e}})$ have same type $n$.
\par By Theorem \ref{DC1}, the cardinality of minimal generating set of $\overline{C(\Gamma_{4})}$ is $2t-2=n-2$. Using induction and from Theorem \ref{TP}, $\beta_{1}(\overline{C(\Gamma_{e}}))= \beta_{1}(\overline{C(\Gamma_{e-1}}))+1=n-2+e$, which is same for $\beta_{1}(C(\Gamma_{e}))$. \qed

\begin{theorem}
 Let n be any positive integer greater than 2, put $s=kn+1, k \in \mathbb{N}$, define
$$\Gamma_{n+1}=\langle s,s-1,\dots,s-n+1,s-n \rangle$$,
now define $\Gamma_{n+1+i}=\langle d,a\Gamma_{n+i}\rangle$, $ 1 \leq i \leq e-(n+1)$ such that $(a,d)=1$ and $d \in \Gamma_{n+1}$	then $C(\Gamma_{e})$ and $\overline{C(\Gamma_{e}})$ have same Cohen Macaulay type $n$ and their first Betti number is $\frac{n(n+4)-2}{2}+(n+1-e)$.															
\end{theorem}

\proof: The semigroup ring $\overline{C(\Gamma_{4}})$ is arithmetically Cohen-Macaulay and its Cohen Macaulay type is $n$ ( See \cite{CN}).  By Theorem \ref{TP}, type of $\overline{C(\Gamma_{5}})$  is $n$. For $e>n+1$, inductively $\Gamma_{e}$ is a simple gluing of $\Gamma_{e-1}$ and again by using Theorem \ref{TP} type of $\overline{C(\Gamma_{e}})$  is $n$. In \cite{CN}, It is given that the type of 
semigroup ring $C(\Gamma_4)$ is $2t-4=n$. Hence $C(\Gamma_{e})$ and $\overline{C(\Gamma_{e}})$ have same type $n$.
\par Since the minimal generator of $\Gamma_{n+1}$ are in arithmetic progression, by remark \ref{CA1}, the first Betti number of $\overline{C(\Gamma_{n+1})}$ is $\frac{n(n+4)-2}{2}$. Using induction and from Theorem \ref{TP}, $\beta_{1}(\overline{C(\Gamma_{e}}))= \beta_{1}(\overline{C(\Gamma_{e-1}}))+1=\frac{n(n+4)-2}{2}+(n+1-e)$, which is same for $\beta_{1}(C(\Gamma_{e}))$. \qed

\begin{theorem}
Let n be an positive integer greater than 2, choose $b$ such that $(n+1,b)=1$ 
$$\Gamma_{n+1}=\langle n+1,n+1+b,n+1+2b,\dots, (n+1)+nb \rangle$$
now define $\Gamma_{n+1+1}=\langle d,a\Gamma_{n+1} \rangle$ such that $(a,d)=1$ and $d \in \Gamma_{n+1}$ ,inductively define  $\Gamma_{e}=\langle d,a\Gamma_{e-1} \rangle, e>n+1$, then $C(\Gamma_{e})$ and $\overline{C(\Gamma_{e}})$ have same Cohen Macaulay type $n$, for $e \geq n+1$ and their first Betti number is $\frac{n(n+1)}{2}+(n+1-e)$.
\end{theorem}											

\proof: Since $\Gamma_{n+1}$ is minimally generated by arithmetic sequence and it is maximal embedding dimension numerical semigroup, hence type of $C(\Gamma)$ is $n$.  Since $\overline{C(\Gamma_{n+1})}$ is associated to arithmetic sequence, it is arithmetically Cohen-Macaulay(see \cite{ACM}). In the proof of Theorem 2.13 in \cite{Bermejo}, it is showed that type of $\overline{C(\Gamma_{n+1}})$ is $n$. For $e>n+1$, inductively $\Gamma_{e}$ is a simple gluing of $\Gamma_{e-1}$ and using Theorem \ref{TP} type of $\overline{C(\Gamma_{e}})$  is $n$. Hence $C(\Gamma_{e})$ and $\overline{C(\Gamma_{e}})$ have same Cohen Macaulay type $n$.	

\par By Remark \ref{CA1}, the first Betti number of $\Gamma_{n+1}$ is $\frac{n(n+1)}{2}$. Using induction and from Theorem \ref{TP}, $\beta_{1}(\overline{C(\Gamma_{e}}))= \beta_{1}(\overline{C(\Gamma_{e-1}}))+1=\frac{n(n+1)}{2}+(n+1-e)$, which is same for $\beta_{1}(C(\Gamma_{e}))$.
\qed


\begin{thebibliography}{00}

\bibitem{SMR}

 P. Gimenez, H. Srinivasan,The structure of the minimal free resolution of semigroup rings
obtained by gluing,Journal of Pure and Applied Algebra
Volume 223, Issue 4, April 2019, Pages 1411-1426


\bibitem{CMC}

J\"{u}rgen-Herzog, Dumitru I. Stamate. "Cohen-Macaulay Criterion for projective monomial Curves via Gr\"{o}bner Bases", Acta Mathematica Vietnamica,2019.

\bibitem{CN} Maria Pia Cavaliere, Gianfranco Niesi."On monomial curves and Cohen-Macaulay type", Manuscripta Mathematica,1983.

\bibitem{BC}
G.-M. Greuel, G. Pfister, A Singular Introduction to Commutative Algebra, Springer-Verlag,
2002.

\bibitem{GHM}
F. Arslan, P. Mete, and M. Sahin. “Gluing and Hilbert functions of monomial curves”. Proc.
Amer. Math. Soc. 137.7 (2009), pp. 2225–2232.

\bibitem{Macaulay}
Daniel R. Grayson and Michael E. Stillman. Macaulay 2, a software system for research in
algebraic geometry. Available at http://www.math.uiuc.edu/Macaulay2

\bibitem{Kunz}
E. Kunz.  “The value-semigroup of a one-dimensional Gorenstein ring”.Proc. Amer. Math.Soc.25(1970), pp. 748–751

\bibitem{RS}
 García-Sánchez, J.C. Rosales, P.A. (2009). Numerical semigroups (First. ed.). New York: Springer. p. 7. ISBN 978-1-4419-0160-6.

\bibitem{JCR}
J.C.Rosales, On presentations of subsemigroups of $\mathbb{N}^{n}$, Semigroup Forum 55(1997) 152-159.

\bibitem{ST}
Richard P. Stanley, Hilbert functions of graded algebras, Adv. Math. 28
(1978), 57-83 [103].

\bibitem{ACM}
D. P. Patil and L. G. Roberts. “Hilbert functions of monomial curves”. J. Pure Appl. Algebra
183.1-3 (2003), pp. 275–292.


%\bibitem{IBS}
%Ping Li, D.P.Patil, Leslie G Roberts, Bases and Ideal Generators for Projective Monomial Curves, Communications in Algebra. 40(2012), 173-191.
\bibitem{STA}
D. I. Stamate. “Betti numbers for numerical semigroup rings”. Multigraded algebra and
applications. Vol. 238. Springer Proc. Math. Stat. Springer, Cham, 2018, pp. 133–157.  

\bibitem{g}
{W. Gastinger, \emph{Über die Verschwindungsideale monomialer Kurven}, PhD thesis, Univ. Regensburg, Landshut
(1989).
}

\bibitem{TOA}
J. Saha, I. Sengupta, and P. Srivastava. “Projective closures of affine monomial curves”. arXiv preprint arXiv:2101.12440, 2021




\bibitem{RS}
 García-Sánchez, J.C. Rosales, P.A. (2009). Numerical semigroups (First. ed.). New York: Springer. p. 7. ISBN 978-1-4419-0160-6.




\bibitem{TA} P. Gimenez, I. Sengupta, H. Srinivasan.
Minimal graded free resolutions for monomial curves defined by arithmetic sequences
J. Algebra, 388 (2013), pp. 294-310.



\bibitem{SG} Jürgen Herzog. Dumitru I. Stamate. "Quadratic numerical semigroups and the Koszul property." Kyoto J. Math. 57 (3) 585 - 612, September 2017. https://doi.org/10.1215/21562261-2017-0007.

\bibitem{BAS} P. Gimenez, I. Sengupta, and H. Srinivasan. “Minimal graded free resolutions for monomial
curves defined by arithmetic sequences”. J. Algebra 388 (2013), pp. 294–310.

\bibitem{Bermejo}
I. Bermejo, E. García-Llorente, I. García-Marco
Algebraic invariants of projective monomial curves associated to generalized arithmetic sequences
J. Symbolic Comput., 81 (2017), pp. 1-19

\bibitem{Watanabe}
 K. Watanabe, Some examples of one dimensional Gorenstein domains, Nagoya Math. J. 49
(1973), 101–109.


\bibitem{Peeva}
I. Peeva, Graded syzygies, in Algebra and applications, vol. 14, Springer,
London, 2011

\bibitem{Rotman}
 J. Rotman. An Introduction to Homological Algebra (Academic Press, 1979

\end{thebibliography}
\end{document}